\documentstyle{article}

\makeatletter

\renewcommand{\theequation}{%
\arabic{section}.\arabic{equation}}
\@addtoreset{equation}{section}

\makeatother

\newtheorem{lemm}{Lemma}[section]
\newtheorem{lem}{Lemma}
\newtheorem{pro}{Proposition}
\newtheorem{thm}{Theorem}
\newtheorem{cor}{Corollary}
\newtheorem{ass}{Assumption}

\def\Z{{\bf Z}}\def\Q{{\bf Q}}
\def\half{{1\over 2}}\def\proof{{\it Proof}}
\def\End{{\rm End}}\def\Ker{{\rm Ker}}\def\Im{{\rm Im}\,}
\def\n{\nonumber}
\def\R{{\cal R}}\def\F{{\cal F}}\def\B{{\cal B}}\def\V{{\cal V}}
\def\N{{\cal N}}\def\W{{\cal W}}\def\H{{\cal H}}
\def\I{{\rm I}}\def\II{{\rm II}}
\def\ka{\kappa_{ij}}
\def\hU{{\hat U}}\def\U{{\cal U}}\def\hUU{{\hat \U}}
\def\e{\epsilon}\def\ba{{\bar a}}\def\ha{{\hat a}}
\def\m{\mu_0}\def\p{\partial}
\def\tE{{\tilde E}}\def\tF{{\tilde F}}\def\tphi{{\tilde \phi}}
\def\hotimes{{\hat \otimes}}
\def\bDelta{{\bar \Delta}}
\def\bPhi{{\bar \Phi}}\def\tPhi{{\tilde \Phi}}\def\tPsi{{\tilde \Psi}}
\def\cPhi{\check{\Phi}}
\def\hiota{{\hat \iota}}\def\tiota{{\tilde \iota}}
\def\trho{{\tilde\rho}}\def\hrho{{\hat \rho}}
\def\hV{{\hat \V}}\def\Y{{\hat Y}}
\begin{document}
\title{Toroidal and level 0 $U'_q(\widehat{sl_{n+1}})$ actions 
on $U_q(\widehat{gl_{n+1}})$ modules}
\author{Kei Miki\\
Department of Mathematics,\\
Graduate School of Science, Osaka University\\
Toyonaka 560, Japan}
\date{}
\maketitle
\begin{abstract}

(1)
Utilizing  a Braid group action on a completion of $U_q(\widehat{sl_{n+1}})$,
an algebra homomorphism from the toroidal algebra  $U_q(sl_{n+1,tor})$
$ (n\ge 2)$ with fixed parameter
to a completion of $U_q(\widehat{gl_{n+1}})$ is obtained.
(2) The toroidal actions by Saito induces 
a level 0 $U'_q(\widehat{sl_{n+1}})$ action
on level 1 integrable highest weight  modules of 
$U_q(\widehat{sl_{n+1}})$. Another 
level 0 $U'_q(\widehat{sl_{n+1}})$ action is defined
 by Jimbo,  et al., in the case $n=1$.
Using   the fact that
the intertwiners of $U_q(\widehat{sl_{n+1}})$ modules are intertwiners of 
toroidal modules for an appropriate comultiplication, 
the relation between  these two level 0 $U'_q(\widehat{sl_{n+1}})$ actions
is clarified.
\end{abstract}

\section*{I. \hskip5mm Introduction}
In \cite{GKV} and \cite{VV1},
a quantum toroidal algebra $U_q(sl_{n+1,tor})$ is introdued.
Up to now several results are obtained on this algebra.
In \cite{VV1}
the connection between toroidal modules and an extension of the double 
affine Hecke algebra \cite{C} is noticed and 
the Schur type duality is obtained.
The vertex representations are  constructed  on level 1
$U_q(\widehat{gl_{n+1}})$ modules by Saito \cite{Saito}.  
In \cite{UT},  toroidal actions are shown to be defined 
on any integrable highest weight module of $U_q(\widehat{gl_{n+1}})$,
using the level - rank duality.
Since the toroidal algebra has  homomorphic images of 
$U_q(\widehat{sl_{n+1}})$ and $U'_q(\widehat{sl_{n+1}})$,
$U_q(\widehat{sl_{n+1}})$ and $U'_q(\widehat{sl_{n+1}})$ 
actions are defined on toroidal modules.
The known  $U'_q(\widehat{sl_{n+1}})$ actions obtained in this way
have  level 0. (See \cite{CP}.) Therefore
level 0 $U'_q(\widehat{sl_{n+1}})$  actions 
on $U_q(\widehat{sl_{n+1}})$ modules are closely related to  toroidal modules.
In  \cite{VV2}\cite{UTS}, the level 1 $U_q(\widehat{sl_{n+1}})$ action
on the fermionic Fock space \cite{KMS}
and the level 0 $U'_q(\widehat{sl_{n+1}})$ action via the affine Hecke algebra
are shown to be combined into a toroidal action.
In \cite{U_0}, motivated by \cite{BPS}, 
a level 0 $U'_q(\widehat{sl_2})$ action is defined
on level 1 integrable $U_q(\widehat{sl_2})$ modules,
utilizing the intertwiners and the representation of the affine Hecke algebra.

In this paper we obtain two results on these problems.
In Sect. III,
utilizing  a Braid group action on a completion of $U_q(\widehat{sl_{n+1}})$,
an algebra homomorphism from the toroidal algebra
$U_q(sl_{n+1,tor})$ $(n\ge 2)$  with fixed parameter
to  a completion of $U_q(\widehat{gl_{n+1}})$ is constructed.
(The paramter depends on the central element.)
 This implies that
any highest weight module of $U_q(\widehat{gl_{n+1}})$ is a toroidal module.
This result corresponds to the fact that the algebra homomorphism
from  $U_q(\widehat{sl_{n+1}})$ to $U_q(gl_{n+1})$ by Jimbo \cite{J}
is neatly expressed in terms of the Braid group action by Lusztig \cite{L}.
In sections IV and V, assuming a triangular decomposition of the
toroidal algebra, 
we consider the level 1 toroidal modules by Saito.
Utilizing   the fact that 
the intertwiners of $U_q(\widehat{sl_{n+1}})$ modules are intertwiners of 
toroidal modules for an appropriate comultiplication, 
we clarify the relation between  the level 0 $U'_q(\widehat{sl_{n+1}})$ action
induced by the toroidal action and that in \cite{U_0}.
Note  that in \cite{BPS} firstly  a Yangian action on  level 1 integrable
highest weight modules of $\widehat{sl_2}$  
is  constructed  in terms of the currents and then
the intertwining property of the vertex operators  is used.
Therefore our approach  is closer to  the original one.
Clarifying  the connection between our results and
\cite{UT}\cite{VV2}\cite{UTS} would be interesting.

Near the  completion of  this work, we found that 
the new Braid group action (Proposition \ref{pro:br}) is obtained in 
\cite{DK}, where it is used in a different manner.

\section*{II. \hskip5mm Definition of algebras}
\addtocounter{section}{1}
\setcounter{subsection}{0}
\setcounter{equation}{0}\addtocounter{section}{1}
\setcounter{subsection}{0}
\setcounter{equation}{0}

Let $q$ be an indeterminate and set $F=\Q(q)$.
Fix $n\ge 2$.

Let $(a_{ij})_{1\le i,j\le n}$ be the Cartan matrix of type $A_n$
and set $\ka=1$.
Let $U_q(\widehat{sl_{n+1}})$  \cite{D}  be the  $F$ algebra defined by the 
generators 
$E_{i,m}$, $F_{i,m}$, $h_{i,r}$, $k_i^{\pm 1}$, $C^{\pm 1}$, $D^{\pm 1}$
$(1\le i\le n,  m\in \Z, r\in \Z\setminus\{0\})$
and the relations,
\begin{eqnarray}
&&C^{\pm 1} \mbox{central},\quad C^{\pm 1} C^{\mp 1}=D^{\pm 1} D^{\mp 1}=1,
 \quad DX_i(z)D^{-1}=X_i(z/q)\quad (X=E,F),\label{eq:top}\\
&&k_i^{\pm 1}k_i^{\mp 1}=1,\quad  D\phi^\pm_i(z)D^{-1}=\phi^\pm_i(z/q),
\label{eq:g0}\\
&&\phi^{\pm}_i(z)\phi^{\pm}_j(w)=\phi^{\pm}_j(w)\phi^{\pm}_i(z),\label{eq:g1}\\
&&{1-q^{-a_{ij}}C^{-1}\ka z/w\over 1-q^{a_{ij}}C^{-1}\ka z/w}
\phi^+_i(z)\phi^-_j(w)=
{1-q^{-a_{ij}}C\ka z/w\over 1-q^{a_{ij}}C\ka z/w}
\phi^-_j(w)\phi^+_i(z),\label{eq:g2}\\
&&\phi^\pm_i(z)E_j(w)\phi^\pm_i(z)^{-1}=
q^{\mp a_{ij}}{1-q^{\pm a_{ij}}C^{-\half\mp\half}(\ka z/w)^{\pm1}
\over 1-q^{\mp a_{ij}}C^{-\half\mp\half}(\ka z/w)^{\pm1}}E_j(w),
\label{eq:g3}\\
&&\phi^\pm_i(z)F_j(w)\phi^\pm_i(z)^{-1}=
q^{\pm a_{ij}}{1-q^{\mp a_{ij}}C^{\half\mp\half}(\ka z/w)^{\pm1}
\over 1-q^{\pm a_{ij}}C^{\half\mp\half}(\ka z/w)^{\pm1}}F_j(w),
\label{eq:g4}\\
&&[E_i(z),F_j(w)]=
{ \delta_{ij}\over q-q^{-1}}
\left(\delta(Cw/z)\phi^-_i(z)-\delta(Cz/w)\phi^+_i(w)\right),\label{eq:ef}\\
&&q^{a_{ij}}(1-q^{-a_{ij}}\ka z/w)E_i(z)E_j(w)=
(1-q^{a_{ij}}\ka z/w)E_j(w)E_i(z), \\
&&q^{-a_{ij}}(1-q^{a_{ij}}\ka z/w)F_i(z)F_j(w)=
(1-q^{-a_{ij}}\ka z/w)F_j(w)F_i(z),\\
&&\mbox{For }i,j \mbox{ such that } a_{ij}=-1\nonumber\\
&&E_i(z_1)E_i(z_2)E_j(w)-(q+q^{-1})E_i(z_1)E_j(w)E_i(z_2)+
E_j(w)E_i(z_1)E_i(z_2)+(z_1\leftrightarrow z_2)=0,\n\\
&&\label{eq:eee}\\
&&F_i(z_1)F_i(z_2)F_j(w)-(q+q^{-1})F_i(z_1)F_j(w)F_i(z_2)+
F_j(w)F_i(z_1)F_i(z_2)+(z_1\leftrightarrow z_2)=0,\n\\
&&\\
&&\mbox{For }i,j \mbox{ such that } a_{ij}=0\nonumber\\
&&[E_i(z), E_j(w)]=0,\quad [F_i(z), F_j(w)]=0,\label{eq:bottom}
\end{eqnarray}
where 	
\begin{eqnarray}
&&E_i(z)=\sum_{m\in\Z}E_{i,m}/z^m, \quad F_i(z)=\sum_{m\in\Z}F_{i,m}/z^m,\n\\
&&\phi^\pm_i(z)=k^{\mp 1}_i\exp\left(\mp(q-q^{-1})
\sum_{r>0} h_{i,\mp r} z^{\pm  r}\right).
\end{eqnarray}
As is well known,
this algebra is also described by the Chevally generators $e_i$,
$f_i$, $k_i^{\pm 1}$ $(0\le i\le n)$ and $D^{\pm1}$. 
Later we need its comultiplication $\Delta_0$ determined by
\begin{eqnarray}
&&\Delta_0(e_i)=e_i\otimes 1+k_i\otimes e_i,\quad 
\Delta_0(f_i)=f_i\otimes k_i^{-1}+ 1\otimes f_i,\n\\
&&\Delta_0(k_i)=k_i\otimes k_i, \quad \Delta_0(D)=D\otimes D.
\end{eqnarray}

$U_q(\widehat{gl_{n+1}})$ \cite{DF}
is defined  to be  the  $F$ algebra generated by 
$E_{i,m}$, $F_{i,m}$, $a_{k,r}$, $t_k^{\pm 1}$, $C^{\pm 1}$, $D^{\pm 1}$
$(1\le i\le n,  1\le k \le n+1, m\in \Z, r\in \Z\setminus\{0\})$
with the defining  relations (\ref{eq:top}), (\ref{eq:ef}--\ref{eq:bottom})
and the following.
\begin{eqnarray}
&&t_k^{\pm 1}t_k^{\mp 1}=1,\quad  D\psi^\pm_k(z)D^{-1}=\psi^\pm_k(z/q),\\
&&\psi^{\pm}_k(z)\psi^{\pm}_l(w)=\psi^{\pm}_l(w)\psi^{\pm}_k(z),\\
&&
{1-q^{-2}C^{-1}z/w\over 1-C^{-1}z/w}
{1-q^{2\theta(k>l)}C^{-1}z/w\over 1-q^{-2\theta(k<l)}C^{-1}z/w}
\psi^+_k(z)\psi^-_l(w)\nonumber \\
&&=
{1-Cz/w\over 1-q^2Cz/w}
{1-q^{2\theta(k>l)}Cz/w\over 1-q^{-2\theta(k<l)}Cz/w}
\psi^-_l(w)\psi^+_k(z),\\
&&\psi^\pm_k(z)E_j(w)\psi^\pm_k(z)^{-1}=
q^{\mp b_{kj}}{1-q^{\pm(k-\half +{1\over 2}b_{kj})}
C^{-\half\mp\half}(z/w)^{\pm1}\over 1-q^{\pm(k-\half-{3\over 2}b_{kj})}
C^{-\half\mp\half}(z/w)^{\pm1}}E_j(w),\n\\
&&\\
&&\psi^\pm_k(z)F_j(w)\psi^\pm_k(z)^{-1}=
q^{\pm  b_{kj}}{1-q^{\pm(k-\half -{3\over 2}b_{kj})}
C^{\half\mp\half}(z/w)^{\pm1}
\over 1-q^{\pm(k-\half+{1\over 2}b_{kj})}
C^{\half\mp\half}(z/w)^{\pm1}}F_j(w).\n\\
&&
\end{eqnarray}
Here $b_{kj}=\delta_{kj}-\delta_{k,j+1}$ and 
$\theta(\cdot)$ is a step function.
In eq. (\ref{eq:ef}), $\phi^\pm_i(z)$ should be understood as follows,
\begin{equation}
\phi_i^\pm(q^iz)=\psi^\pm_i(z)/\psi^\pm_{i+1}(z).
\end{equation}

Let  $\mu\in F^\times$
and  $(a_{ij})_{0\le i,j\le n}$ be the Cartan matrix of type $A^{(1)}_n$
and set $\ka =1$ $((i,j)\ne (n,0), (0,n))$, 
$\kappa_{n0}=\kappa_{0n}^{-1}=\mu$.
Let $U_{q,\mu}(sl_{n+1,tor})$ \cite{VV1} be 
the $F$ algebra defined by the generators
$E_{i,m}$, $F_{i,m}$, $h_{i,r}$, $k_i^{\pm 1}$, $C^{\pm 1}$, $D^{\pm 1}$
$(0\le i\le n,  m\in \Z, r\in \Z\setminus\{0\})$
and the relations (\ref{eq:top}--\ref{eq:bottom}).
Note that we include  only one scaling element $D$ and its inverse among the 
generators.
Set $\m=q^{n+1}$. For $\mu=(\m C^2)^{\pm 1}$, we also define 
$U_{q,(\m C^2)^{\pm 1}}(sl_{n+1,tor})$ similarly.
Later we shall often use $p=\m\mu$ in stead of $\mu$.

Hereafter we shall write $U$ and $\U_\mu$ for $U_q(\widehat{sl_{n+1}})$
and $U_{q,\mu}(sl_{n+1,tor})$, respectively.
We shall identify $U_q(\widehat{sl_{n+1}})$ with the subalgebra of 
$U_q(\widehat{gl_{n+1}})$ generated by $E_i(z)$, $F_i(z)$, $\phi^\pm_i(z)$
$(1\le i\le n)$, $C^{\pm 1}$ and $D^{\pm 1}$.

\section*{III. \hskip5mm  Homomorphism  from $\U_{(\m C^2)^{\pm 1}}$ to a
 completion of $U_q(\widehat{gl_{n+1}})$}
\addtocounter{section}{1}
\setcounter{subsection}{0}
\setcounter{equation}{0}

\subsection{Braid group action on a completion of $U_q(\widehat{sl_{n+1}})$}
For $k,l\in \Z_{\ge 0}$, put
\begin{equation}
U_{kl}=\sum_{r\ge k\atop s\ge l } U_{-r}U U_s,
\end{equation}
where
\begin{equation}
U_r=\{y\in U| DyD^{-1}=q^ry\}, \quad (r\in \Z).
\end{equation}
Then
$U$ is a   topological  algebra which has
$(U_{kl})$ as a    fundamental system of  neighborhoods of the origin
and is separated thanks to the  triangular decomposition of $U$ \cite{Ro}.
Let us denote its completion by $\hU$.

Let $\eta$ be a continuous  algebra anti-automorphism of  
$\hU$ determined by
\begin{eqnarray}
&&E_i(z)\mapsto E_i(z^{-1}),\quad F_i(z)\mapsto F_i(z^{-1}),\quad 
\phi^\pm_i(z)\mapsto \phi^\mp_i(Cz^{-1})\nonumber\\
&&C^{\pm 1}\mapsto C^{\pm 1},\quad D^{\pm 1}\mapsto D^{\pm 1}.
\end{eqnarray}

\begin{pro}\label{pro:br}
(1) For each $i$ ($1\le i\le n$), there exists 
a continuous   algebra automorphism $T_i:\hU\to  \hU$
determined by
\begin{eqnarray}
&&E_i(z)\mapsto -F_i(z/Cq^2)\phi^-_i(z/q^2)^{-1},\quad 
F_i(z)\mapsto -\phi^+_i(z/q^2)^{-1}E_i(z/Cq^2),\nonumber\\
&&\phi^\pm_i(z)\mapsto \phi^\pm_i(z/q^2)^{-1},\quad 
C^{\pm 1}\mapsto C^{\pm 1}, \quad D^{\pm 1}\mapsto D^{\pm 1},\\
&&\nonumber\\
&&E_j(z)\mapsto \oint{du\over u}\left(E_i(u)E_j(z)-
q{1-q^{-1}u/z\over 1-qu/z}E_j(z)E_i(u)\right),\nonumber\\
&&F_j(z)\mapsto \oint{du\over u}\left(F_j(z)F_i(u)-
q^{-1}{1-qz/u\over 1-q^{-1}z/u}F_i(u)F_j(z)\right),\nonumber\\
&&\phi^\pm_j(z)\mapsto \phi^\pm_j(z)\phi^\pm_i(z/q)\quad 
\mbox{ when } a_{ij}=-1,\\
&&\nonumber\\
&&E_j(z)\mapsto E_j(z),\quad 
F_j(z)\mapsto F_j(z),\quad \phi^\pm_j(z)\mapsto \phi^\pm_j(z)
\quad \mbox{ when } a_{ij}=0.
\end{eqnarray}
Here $\displaystyle \oint {du\over u}$
denotes the operation which picks out the coefficient
of $u^0$.

(2) $T_i^{-1}=\eta\circ T_i\circ \eta.$

(3) $T_i$'s satisfy the following relations.
\begin{eqnarray}
&&T_iT_jT_i=T_jT_iT_j\quad \mbox{when   } a_{ij}=-1,\nonumber\\ 
&&T_iT_j=T_jT_i \quad \mbox{when   } a_{ij}=0.
\end{eqnarray}

(4) Let $W$ be the Weyl group of $sl_{n+1}$ and 
$s_i$ be a  reflection with respect to a simple root $\alpha_i$.
For a reduced expression $w=s_{i_1}\cdots s_{i_m}\in W$, set
$T_w=T_{i_1}\cdots T_{i_m}$. If $w\alpha_k=\alpha_l$,  then 
\begin{equation}
T_w(X_k(z))=X_l(q^{-r/2}z)  \quad ( X=E, F, \phi^\pm).
\end{equation}
Here $r$ is the length of the shortest element(s) of the set
$\{w'\in W\mid w'\alpha_k=\alpha_l\}$.
\end{pro}

\subsection{Homomorphism from $\U_{(\m C^2)^{\pm 1}}$ to 
a completion of $U_q(\widehat{gl_{n+1}})$}
Similarly to $U$, we introduce a separated linear topology
on   $U_q(\widehat{gl_{n+1}})$ and denote its completion by
$\hU_q(\widehat{gl_{n+1}})$.
Utilizing the above Braid group action, we obtain the following theorem.
\begin{thm}\label{thm:hom}
Set
\begin{eqnarray}
&&\tE_0^\e(z)=-T_n^\e\cdots T_1^\e E_1(zq^\e),\quad
 \tF_0^\e(z)=-T_n^\e\cdots T_1^\e F_1(zq^\e),\n\\
&&\tphi^{\pm,\e}_0(z)=T_n^\e\cdots T_1^\e\phi^\pm_1(zq^\e),\quad
(\e=\pm 1).
\end{eqnarray}
For $\e=\pm 1$,
 there exists an  algebra homomorphism $f_\e:\U_{(\m C^2)^\e}
\to \hU_q(\widehat{gl_{n+1}})$  determined by
\begin{eqnarray}
&&E_i(z)\mapsto E_i(z),\quad F_i(z)\mapsto F_i(z),\quad
\phi^\pm_i(z)\mapsto \phi^\pm_i(z),\quad (1\le i\le n),\n\\
&&C\mapsto C,\quad D\mapsto D,\n\\
&&E_0(z)\mapsto \psi^+_{n+1}(z/ \mu_\e)^\e
\tE_0^\e(z)\psi^-_{n+1}(Cz/ \mu_\e)^{-\e},\n\\
&&F_0(z)\mapsto \psi^+_{n+1}(Cz/ \mu_\e)^{-\e}
\tF_0^\e(z)\psi^-_{n+1}(z/ \mu_\e)^\e,\n\\
&&\phi^\pm_0(z) \mapsto\tphi^{\pm,\e}_0(z)
 \psi^{\pm}_{n+1}(z/ C \mu_\e)^\e \psi^\pm_{n+1}(Cz/ \mu_\e)^{-\e},
\end{eqnarray}
where $\mu_\e=\m (\m C)^\e$.
\end{thm}

\proof.
As an example, we shall check the Serre relation (\ref{eq:eee})
with $(i,j)=(0, n)$ for the case $\e=1$.
Set $w=s_n\cdots s_1$. 
Thanks to Proposition \ref{pro:br} (4), we get
\begin{equation}
T_w^{-1}E_{n-1}(z/q)=E_n(z),\quad 
T_w^{-1}E_n(z/q)=
-\tphi_0^{+,+}(\mu z/C)\tE^{+}_0(\mu z)\tphi^{-,+}_0(\mu z)^{-1}
\end{equation}
where $\mu=\m C^2$.
Therefore applying $T_w^{-1}$ to 
the Serre relation (\ref{eq:eee}) with $(i,j)=(n,n-1)$
of $U_q(\widehat{gl_{n+1}})$ and using the relations among $\psi^{\pm}_k(z)$,
$E_i(z)$ and $F_i(z)$, 
we obtain the desired relation.

\begin{cor}
Any highest weight module of $U_q(\widehat{gl_{n+1}})$ on which $C$ acts as 
$C_0\in F^\times$ is a $\U_{(\m C_0^2)^{\pm 1}}$ module.
\end{cor}

\section*{IV. \hskip5mm Intertwiners of toroidal modules}
\addtocounter{section}{1}
\setcounter{subsection}{0}
\setcounter{equation}{0}

\subsection{Level 1 toroidal modules by Saito}
Hereafter we consider only $\U_\mu$ ($\mu\in F^\times$).
Let $V(\Lambda_i)$ $(0\le i\le n)$ be the irreducible
highest weight module of $U$ with highest weight $\Lambda_i$
and $v_{\Lambda_i}$ be its highest weight vector.

Let $\B$ be the $F$ algebra generated by $b_r$
 $(r\in \Z-{0})$ and $D^{\pm 1}$
with the relations
\begin{eqnarray}
&&D^{\pm 1} D^{\mp 1}=1,\quad D b_rD^{-1}=q^{r}b_r,\n\\
&&[b_r,b_s]=\delta_{r+s,0}{[r]^2[(n+2)r]\over r[(n+1)r]},
\end{eqnarray}
where $\displaystyle [m]={q^m-q^{-m}\over q-q^{-1}}$.
Let  $\B_-$ denote its Fock space $\B/I$, where $I$ is the left ideal 
generated by $b_r$ $(r>0)$ and $D-1$.
For $m\equiv j$ mod $n+1$, set $W_m=V(\Lambda_j)\otimes \B_-$.
Let $\alpha_i$'s be  the simple roots of  $\widehat{sl_{n+1}}$
and $(\,\mid\,)$ be the standard 
symmetric bilinear form on the weight space of $\widehat{sl_{n+1}}$
normalized by $(\alpha_i|\alpha_j)=a_{ij}$.
Set
$$
\e_k=\left(\sum_{i=k}^n(n+1-i)\alpha_i-\sum_{i=1}^{k-1}i\alpha_i\right)
/(n+1),\quad (1\le k\le n+1).
$$
For $a\in F^\times$ and $1\le k\le n+1$,
let  $a^{\partial_{\epsilon_k}}$   be a linear operator on $W_m$
such that $a^{\partial_{\epsilon_k}}(v\otimes b) 
=a^{(\epsilon_k|\nu)+{m\over n+1}}v\otimes b$,
where $v$ is  a weight vector with weight $\nu$.
Then   $W_m$  is a $U_q(\widehat{gl_{n+1}})$
 module by the  map
\begin{eqnarray}
&&E_i(z)\mapsto E_i(z)\otimes 1,\quad F_i(z)\mapsto F_i(z)\otimes 1,\n\\
&&\psi^\pm_k(z)\mapsto q^{\mp \partial_{\epsilon_k}}
\exp\left(\mp(q-q^{-1})\sum_{r>0} \ha_{k, \mp r} z^{\pm r}\right),\n\\
&&C\mapsto q,\quad D\mapsto  D\otimes D,
\end{eqnarray}
where $\ha_{k,r}=\ba_{k,r}\otimes 1+1\otimes b_r$ with 
\begin{equation}
\ba_{k,r}={1\over [(n+1)r]}\left(
\sum_{i=k}^n[(n+1-i)r]h_{i,r}-q^{-(n+1)r}\sum_{i=1}^{k-1} [ir]h_{i,r}\right).
\end{equation}

The results by Saito \cite{Saito}
on level 1 toroidal modules  can be stated as follows.
\begin{pro}
Set $a=(-1)^{n-1}q^n$.

(1) For $\e=\pm 1$, the $U$ module structure  on $V(\Lambda_j)$
 $(0\le j\le n)$ is extended to a $\U_{\m^\e}$ module structure
 by 
\begin{equation}
E_0(z)=a \tE^\e_0(z),\quad F_0(z)=a^{-1}\tF_0^\e(z),\quad 
\phi^\pm_0(z)=\tphi^{\pm,\e}_0(z).
\end{equation}
(2) 
Letting $\alpha_r$ $(r\ne 0)$ be the elements of $F^\times$ such that 
\begin{equation}
\alpha_r\alpha_{-r}={\mu^r+\mu^{-r}-\m^r-\m^{-r}\over 
(q^r-q^{-r})(q^r\m^r-q^{-r}\m^{-r})},
\end{equation}
set
\begin{eqnarray}
&&X^\pm(z)=
\exp\left( \mp(q-q^{-1})\sum_{r>0}x_{\mp r}z^{\pm r}\right),\n\\
&&x_r=(q\m^2)^r{(\mu/\m)^r-1\over q^{2r}-1}\ba_{n+1,r}\otimes 1
+\alpha_r 1\otimes b_r.
\end{eqnarray}
Then the $U$ module structure on   $W_m$  is extended to  
a $\U_{\mu}$ module structure by 
\begin{eqnarray}
&& E_0(z)=a X^+(z)(\tE_0^+(z)\otimes 1)
X^-(qz)^{-1}(\m/\mu)^{\p_{\e_{n+1}}},\n\\
&&F_0(z)=a^{-1}(\mu/\m)^{\p_{\e_{n+1}}}X^+(qz)^{-1}
(\tF_0^+(z)\otimes 1)X^-(z),\n\\
&&\phi^\pm_0(z)={X^\pm(q^{-1}z)\over X^\pm(qz)}
(\tphi^{\pm, +}_0(z)\otimes 1).
\end{eqnarray}
\end{pro}
\vskip3mm
We have chosen the above normalization of $E_0(z)$ and $F_0(z)$ for 
a later convenience. Note that in the case $\mu=\m q^2$ the above result
coincides with Theorem \ref{thm:hom}.

\subsection{Level 0 toroidal modules}
Set $V=F^{n+1}$ and let $v_1,\cdots,v_{n+1}$ be its canonical basis.
$V[x,x^{-1}]$ is a $U$ module by the following.
\begin{eqnarray}
&&E_i(q^i z)=\delta(z/x)E_{i\,i+1},\quad
 F_i(q^i z)=\delta(z/x)E_{i+1\,i},\n\\
&&\phi^+_i(q^i z)=\sum_{k\ne i, i+1}E_{kk}+
q^{-1}{1-q^2z/x\over 1-z/x}E_{ii}
+q{1-q^{-2}z/x\over 1-z/x}E_{i+1\,i+1},\n\\
&&\phi^-_i(q^i z)=\sum_{k\ne i, i+1}E_{kk}
+q{1-q^{-2}x/z\over 1-x/z}E_{ii}
+q^{-1}{1-q^2x/z\over 1-x/z}E_{i+1\,i+1},\n\\
&&\hskip10cm(1\le i\le n),\n\\
&&C^{\pm 1}= 1,\quad D^{\pm 1}= q^{\pm \vartheta}.
\end{eqnarray}
Here $E_{ij}$'s  are matrix units,
$\delta(z)=\sum_{m\in \Z} z^m$, $x$ acts as the multiplication,
and for $a\in F^\times$  and $v(x)\in V[x, x^{-1}]$,
$a^\vartheta v(x)=v(ax)$.
We shall denote this representation by $(\pi, V_x)$.

Set $p=\m\mu$.
The $U$ module structure on $V_x$ is extended to a $\U_\mu$
structure \cite{VV1} by
\begin{eqnarray}
&&E_0(z)=a p^\vartheta\delta(z/x)E_{n+1\,1},\quad
F_0(z)=a^{-1}\delta(z/x)E_{1\, n+1}p^{-\vartheta},\n\\
&&\phi^+_0(z)=\sum_{k\ne 1, n+1} E_{kk}+
q^{-1}{1-q^2z/px\over 1-z/px}E_{n+1\, n+1}
+q{1-q^{-2}z/x\over 1-z/x}E_{11},\n\\
&&\phi^-_0(z)=\sum_{k\ne 1, n+1}E_{kk}
+q{1-q^{-2}px/z\over 1-px/z}E_{n+1\, n+1}
+q^{-1}{1-q^2x/z\over 1-x/z}E_{11},\n\\
&&
\end{eqnarray}
where  $a\in F^\times$.
We shall denote this representation by $(\pi_a, V_a)$.

\subsection{ Bialgebra structure of $\U_\mu$}

In this subsection, fixing $\mu\in F^\times$, we omit the subscript $\mu$
of $\U_\mu$. 
Define $\U_r$  $(r\in \Z)$ similarly to $U_r$.  We assume the following.
\begin{ass}
(1) $U$ is identified with the subalgebra of $\U$
generated by $E_i(z)$, $F_i(z)$, $\phi^\pm_i(z)$ $(1\le i\le n)$,
$C^{\pm 1}$ and $D^{\pm 1}$ by the  algebra homomorphism
from $U$ to $\U$  determined by
\begin{equation}
E_i(z)\mapsto E_i(z),\quad
F_i(z)\mapsto F_i(z),\quad
\phi^\pm_i(z)\mapsto \phi^\pm_i(z)\quad
C\mapsto C,\quad D\mapsto D.
\end{equation}

(2) 
$\U$ has subalgebras $\U^-$, $\U^0$ and $\U^+$ such that 
the multiplication map $\U^-\otimes \U^0\otimes \U^+\to \U$
is an isomorphism of  vector spaces and 
\begin{equation}
U^\pm \subset\U^\pm=\oplus _{r\ge 0} \U^\pm\cap \U_{\pm r},\quad 
U^0\subset \U^0\subset \U_0,
\end{equation}
where
$U^+$, $U^-$ and $U^0$ are the subalgebras of $U$ generated by
$e_i$ $(0\le i\le n)$,  $f_i$ $(0\le i\le n)$,
and $k_i^{\pm 1}$ $(0\le i\le n)$  and $D^{\pm 1}$, respectively.
\end{ass}

For $r\in \Z$ setting
\begin{equation}
\U^{\otimes N}{}r=\{y\in \U^{\otimes N}\mid  D^{\otimes N}
y(D^{\otimes N})^{-1}=q^ry\},
\end{equation}
define  $\U^{\otimes N}{}_{kl}$ similarly  to  $U_{kl}$.
We introduce a linear topology on  $\U^{\otimes N}$ 
by letting $(\U^{\otimes N}{}_{0l})_{l\ge 0}$
be a fundamental system of neighborhoods of the origin. 
(Hereafter we shall simply say  `introduce a linear topology
by $(\U^{\otimes N}{}_{0l})$'.)
Then $\U^{\otimes N}$ is a separated topological  algebra.
We denote its completion by $\hUU$ ($N=1$) and $\U^{\hotimes N}$
$(N\ge 2)$. 
Note that for $N\ge 2$, the topology on $\U^{\otimes N}$
is equivalent to the tensor product topology.

Let $\bDelta_\I$ and $\bDelta_\II$ be the continuous  algebra homomorphism
from $\hUU$ to $\U\hotimes \U$ determined  by
\begin{eqnarray}
&&\bDelta_\I(C)=\bDelta_\II(C)=C\otimes C,\quad
\bDelta_\I(D)= \bDelta_\II(D)=D\otimes D,\n\\
&&\bDelta_\I(\phi^+_i(z))=\phi^+_i(z/C_2)\otimes \phi^+_i(z),\quad
\bDelta_\I(\phi^-_i(z))=\phi^-_i(z)\otimes \phi^-_i(z/C_1),\n\\
&&\bDelta_\I(E_i(z))=E_i(z)\otimes 1+\phi_i^-(z)\otimes E_i(z/C_1),\quad
\bDelta_\I(F_i(z))=1\otimes F_i(z)+F_i(z/C_2)\otimes \phi_i^+(z),\n\\
&&\bDelta_\II(\phi^+_i(z))=\phi^+_i(z/C_2)\otimes \phi^+_i(z),\quad
\bDelta_\II(\phi^-_i(z))=\phi^-_i(z)\otimes \phi^-_i(z/C_1),\n\\
&&\bDelta_\II(E_i(z))=E_i(z)\otimes 1+\phi_i^+(Cz)\otimes E_i(C_1z),\quad
\bDelta_\II(F_i(z))=1\otimes F_i(z)+F_i(C_2z)\otimes \phi_i^-(Cz),\n\\
&&
\end{eqnarray}
where $C_1=C\otimes 1$ and $C_2=1\otimes C$.
Let further $\e:\hUU\to F$ denote the continuous  algebra homomorphism
 determined by
\begin{equation}
\e(E_i(z))=\e(F_i(z))=0,\quad \e(\phi^\pm_i(z))=\e(C)=\e(D)=1.
\end{equation}
Here $F$ is given a discrete topology.
Then $(\hUU,\bDelta_i,\e)$ $(i=\I,\II)$ are bialgebras.
These are straightforward generalizations  of the bialgebra structures found
by Drinfeld for $\hU$.
Let $\R=\R^+\R^0\R^-$ be the Gauss decomposition of the universal
$R$ matrix of $U$ with $\Delta_0$ as a comultiplication \cite{KT}.
Letting
\begin{equation}
\F_\I=\R^-,\quad \F_\II=\sigma(\R^+)^{-1},\quad(\sigma(a\otimes b)=b\otimes a),
\end{equation}
set
\begin{equation}
\Delta_i(\cdot)=\F_i^{-1}\bDelta_i(\cdot)\F_i,\quad
(i=\I,\II).
\end{equation}
It is known \cite{KT} that $\Delta_\I(y)=\Delta_\II(y)=\Delta_0(y)$
for $y\in U$ and that $(\hUU,\Delta_i,\e)$  $(i=\I, \II)$
 are bialgebras.
Set $\Delta^{op}_\I=\sigma\circ\Delta_\I$ and 
$\F_\I^{op}=\sigma(\F_\I)$.
For $(\Delta,\F)=(\Delta^{op}_\I, \F_\I^{op})$,
$(\Delta_\II, \F_\II)$, etc.,
we define $\Delta^{(N)}:\hUU\to \U^{\hotimes N+1}$ $(N\ge 1)$ by
\begin{equation}
\Delta^{(1)}=\Delta,\quad \Delta^{(N)}=(\Delta\hotimes 1^{\hotimes N-1})
\circ \Delta^{(N-1)}\quad (N\ge 2).
\end{equation}
and  $\F^{(N)}\in \U^{\hotimes N}$ $(N\ge 2)$  by
\begin{equation}
\F^{(2)}=\F,\quad 
\F^{(N)}=(\F^{(2)})_{12}(\Delta_\I^{op}\hotimes 1^{\hotimes N-2})\F^{(N-1)},
\quad (N\ge 3).
\end{equation}
Then  the following holds.
\begin{equation}
\Delta^{(N-1)}(y)=
\F^{(N)}{}^{-1}\bDelta^{(N-1)}(y) \F^{(N)},\quad (y\in \U). 
\label{eq:delta}
\end{equation}
Later we need the following property of $\F^{op}_\I$ and $\F_\II$
\cite{KT}. Set $U^{\ge 0}=U^0U^+$ and $U^{\le0}=U^-U^0$.
Set further $\bar{Q}=\oplus_{i=1}^n \Z \alpha_i$,
$\bar{Q}_{+}=-\bar{Q}_-=\oplus_{i=1}^n \Z_{\ge 0} \alpha_i\setminus\{0\}$ 
and $\delta=\sum_{i=0}^n \alpha_i$.
For $m\in \Z$ and $\lambda\in \bar{Q}$, define
$$
U^a_{m\delta+\lambda}=\{y\in U^a\mid
DyD^{-1}=q^m y,\quad k_iyk_i^{-1}=q^{(\alpha_i|\lambda)}y,\quad 
(1\le i\le n)\},
$$
where $a=\ge 0$, $\le 0$. Let $M_\pm$ be the closure 
of $\sum_{m\ge 0}\sum_{\lambda\in\bar{Q}_\pm}
U^{\le 0}_{-(m\delta+\lambda)}\otimes U^{\ge 0}_{m\delta+\lambda}$.
Then
\begin{equation}
\F_\I^{op}-1\in M_-,\quad \F_\II-1\in M_+.\label{eq:wt}
\end{equation}

For $N\ge 2$ and  $m\ge 0$, set
\begin{equation}
\U^{\otimes N}{}[m]=\sum_{\sum_{i=2}^N (i-1) m_i\ge m}
\U\otimes \U_{m_2}\otimes\cdots\otimes \U_{m_N}.
\end{equation}
Let $\widetilde{\U^{\otimes N}}$ $(N\ge 2)$ be the $F$ algebra
$\U^{\otimes N}$ on which a linear topology is introduced by
$(\U^{\otimes N}[m])$.  Then 
$\widetilde{\U^{\otimes N}}$ is a separated topological  algebra. We denote its
completion by $\widehat{\widetilde{\U^{\otimes N}}}$.

\begin{lem}\label{lem:U}
(1) Let $\tau$ be the identity map from $\widetilde{\U^{\otimes N}}$ to 
$\U^{\otimes N}$. Let further
 ${\hat \tau}:\widehat{\widetilde{\U^{\otimes N}}} \to \U^{\hotimes N}$ 
be the continuous extension of the continuous  algebra homomorphism $\tau$.
Then ${\hat \tau}$ is injective.

(2) If we identify $\,\widehat{\widetilde{\U^{\otimes N}}}$ 
with a subalgebra of  $\,\U^{\hotimes N}$ 
via the map ${\hat \tau}$, then the following holds.
\begin{equation}
\Delta_\I^{op (N-1)}(\U)\subset \widehat{\widetilde{\U^{\otimes N}}},\quad
\Delta_\II^{(N-1)}(\U)\subset \widehat{\widetilde{\U^{\otimes N}}}.
\end{equation}
\end{lem}

\proof.
(1)
For $l\in \Z_{\ge 0}$, set
\begin{equation}
M_l=\sum_{\sum l_i=l}\U^{-\otimes N}\otimes \U^{0\otimes N}
\otimes \U^{+}_{l_1}\otimes \cdots\otimes \U^{+}_{l_N}.
\end{equation}
Under the identification  of $\U^{\otimes N}$ with $\U^{-\otimes N}\otimes
 \U^{0\otimes N}\otimes \U^{+\otimes N}$, the following holds.
\begin{equation}
\U^{\otimes N}[m]=\oplus_l\left(\U^{\otimes N}[m]\cap M_l \right),\quad
\U^{\otimes N}_{0m}=\oplus_{l\ge m}M_l.
\end{equation}
From this, we obtain 
\begin{equation}
\cap_{0l}\left(\U^{\otimes N}[m]+\U^{\otimes N}_{0l}\right)=
\U^{\otimes N}[m].
\end{equation}
Utilizing the last equality,
it is easy to show that any Cauchy sequence in $\widetilde{\U^{\otimes N}}$
which converges to $0$ in $\U^{\hotimes N}$
converges to $0$ in $\widehat{{\widetilde{\U^{\otimes N}}}}$.

(2) The claim can be easily checked for $\bDelta_\I^{op}$ and $\bDelta_\II$.
Since  thanks to (\ref{eq:wt})
$\F^{op (N)}_\I$, $\F^{(N)}_\II\in \widehat{\widetilde{\U^{\otimes N}}}$ ,
 we obtain the claim.

\vskip3mm

Let $\W_i$ $(1\le i\le N)$ be $\U$ modules such that 
\begin{equation}
\W_i=\oplus_{m\in \Z} \W_{i,m},\quad \W_{i,m}=\{w\in \W_i\mid Dw=q^m w\},
\end{equation}
and set $\W=\W_1\otimes \cdots \otimes \W_N$. For $m\in \Z$ setting
\begin{equation}
\W [m]= \sum_{\sum_{i=2}^N (i-1)  m_i \ge m}
\W_1\otimes \W_{2, m_2}\otimes \cdots \otimes \W_{N, m_N},
\end{equation}
introduce a separated linear topology on $\W$ by $(\W[m])$.
Then ${\hat \W}$, the completion of $\W$, is a topological 
$\widehat{\widetilde{\U^{\otimes N}}}$ module.  Therefore,
thanks to Lemma \ref{lem:U}, ${\hat \W}$
is a  $\U$ module via the comultiplication
$\Delta^{op}_\I$ or $\Delta_\II$. 
Note that $\hat{\W}=\W$
in the case $N=2$ and $\W_{N, m}=\{0\}$ $(m\gg 0)$.

\subsection{Intertwiners of toroidal modules}

Let $\Phi^{(j)}$ and  $\Psi^{(j)}$ ($0\le j\le n$) denote the intertwiners of 
$U$ modules  determined by
\begin{eqnarray}
&&\Phi^{(j)}:V(\Lambda_j)\otimes V_x \to V(\Lambda_{j+1}),\quad 
\Psi^{(j)}: V_x\otimes V(\Lambda_j)\to V(\Lambda_{j+1}),\n\\
&&\Phi^{(j)}(v_{\Lambda_j}\otimes v_{j+1})=v_{\Lambda_{j+1}},\quad
\Psi^{(j)}(v_{j+1}\otimes v_{\Lambda_j})=v_{\Lambda_{j+1}},
\end{eqnarray}
where 
$V(\Lambda_{j})\otimes V_x$ and $V_x \otimes V(\Lambda_{j})$ 
are $U$ modules via the comultiplication $\Delta_0$
and $\Lambda_{n+1}$ should be understood as $\Lambda_0$.

As usual we define their components by 
\begin{equation}
\Phi^{(j)}_{k,m}u=\Phi^{(j)}(u\otimes v_k x^m),\quad
 \Psi^{(j)}_{k,m}u=\Psi^{(j)}(v_k x^m\otimes u),
\end{equation}
and  set 
\begin{equation}
\Phi^{(j)}_k(z)=\sum_{m\in \Z}\Phi^{(j)}_{k,m}z^{-m},\quad
\Psi^{(j)}_k(z)=\sum_{m\in \Z}\Psi^{(j)}_{k,m}z^{-m}.
\end{equation}

\begin{pro}\label{pro:vo}
 For $m\in \Z$, set $a_m=q^{2m}$. 
In the following,   give 
$V(\Lambda_j)\otimes V_a$  and $W_m\otimes V_a$
(resp.  $V_a\otimes V(\Lambda_j)$ and $V_a\otimes W_m$)
 toroidal module structures  via the comultiplication 
$\Delta_\I$ (resp. $\Delta_\II$).

(1)  For $0\le j\le n$, 
$\Phi^{(j)}:V(\Lambda_j)\otimes V_{a_j}\to V(\Lambda_{j+1})$  and
$\Psi^{(j)}:V_{a_j}\otimes V(\Lambda_j)\to V(\Lambda_{j+1})$ 
are intertwiners of $\U_{\m}$ modules.

(2) Let $p=\m\mu\ne 1$.
Setting 
\begin{equation}
Y^\pm(z)=\exp\left(\mp (q-q^{-1})
 \sum_{r>0}{p^{r\pm r\over 2}\over 1-p^r}
\alpha_{\mp r} b_{\mp r} z^{\pm r}\right),
\end{equation}
put
\begin{equation}
\Xi(z)=Y^+(qz)Y^-(z),\quad \Sigma(z)=Y^+(z) Y^-(qz).
\end{equation}
For $m=(n+1)s +j$ ($s\in \Z, 0\le j\le n$), define 
$\tPhi^{(m)}:W_m\otimes V_{a_m/p^s}\to W_{m+1}$ and 
$\tPsi^{(m)}:V_{a_m/p^s} \otimes W_m\to W_{m+1}$
by the  generating series
\begin{equation}
\tPhi^{(m)}_k(z)= \Phi^{(j)}_k(z)\otimes \Xi(z) \mbox{  and  }
\tPsi^{(m)}_k(z)= \Psi^{(j)}_k(z)\otimes \Sigma(z),
\end{equation}
respectively. Here
$\tPhi^{(m)}_k(z)$ and $\tPsi^{(m)}_k(z)$ are defined 
similarly to  $\Phi^{(j)}_k(z)$ and $\Psi^{(j)}_k(z)$.
Then 
$\tPhi^{(m)}$ and  $\tPsi^{(m)}$  are intertwiners  of $\U_\mu$ modules.
\end{pro}

\proof.
We consider the case $\tPhi^{(m)}$. 
Give $W_m\otimes V_a$ (resp. $V(\Lambda_j)\otimes V_x$)
a $\U_\mu$ (resp. $U$) module structure via the comultiplication
$\bDelta_\I$.
Set $\cPhi^{(m)}=\tPhi^{(m)}\F^{-1}_\I$ and 
$\bPhi^{(j)}=\Phi^{(j)} \F_\I^{-1}$.
Then $\cPhi^{(m)}_k(z)=\bPhi^{(j)}_k(z)\otimes \Xi(z)$.
Utilizing 
the explicit expression of the intertwiner 
$\bPhi^{(j)}:V(\Lambda_j)\otimes V_x\to V(\Lambda_{j+1})$ in terms of bosons 
\cite{DI}, it is straightforward
to show  $\cPhi^{(m)}:W_m\otimes V_{a_m/p^s}\to
W_{m+1}$ is an intertwiner of $\,\U_\mu$ modules.
The claim follows from this.

\section*{V. \hskip5mm  Two level 0  $U_q'(\widehat{sl_{n+1}})$ actions}
\addtocounter{section}{1}
\setcounter{subsection}{0}
\setcounter{equation}{0}

\subsection{ Completions of $V_x^{\otimes N}$}

Following  \cite{U_0}, we introduce two completions of 
$V_x^{\otimes N}$ and the maps from them to $W_N$.

For $m\in \Z$, let $\V_N[m]$ signify $V_x^{\otimes N}[m]$ 
in subsection 4.3 and set 
\begin{eqnarray}
&&\V_N[[m]]={\rm span}\{v_{\e_1}x^{m_1}\otimes \cdots v_{\e_N}x^{m_N}\mid
{\rm max}\{m_1,\cdots, m_N\}\ge m\}.
\end{eqnarray}
Introduce  two linear topologies  on the vector space $V_x^{\otimes N}$ 
by  $(\V_N[m])$   and  by  $(\V_N[[m]])$.
We denote the thus obtained separated topological vector spaces by $\V'_N$
and $\V_N$, respectively,  and
let $\hV'_N$ and  $\hV_N$ signify their  completions.
Let $\iota_N$ be the identity map from $\V'_N$ to $\V_N$.
Let further $\hiota_N:\hV'_N\to \hV_N$
denote the continuous extension of the continuous linear map $\iota_N$.
(This map is shown to be injetive as in Lemma \ref{lem:U}.)
Note that our completions are a little bit different from those in \cite{U_0}.

Set 
\begin{equation}
v_{\e_1,\cdots,\e_N}(z_1,\cdots,z_N)=\sum_{m_i} 
v_{\e_1}x^{m_1}\otimes \cdots \otimes v_{\e_N}x^{m_N}
z_1^{-m_1}\cdots z_N^{-m_N}.
\end{equation}
Let $\N'_N$ (resp. $\N_N$) be the closure  in $\hV'_N$  (resp. $\hV_N$)
of the span of the coefficients of the following  generating series, 
\begin{eqnarray}
&&v_{\cdots ,\e_i,\e_{i+1},\cdots}(\cdots,z_i,z_{i+1},\cdots)
-(1-q^2){(z_i/z_{i+1})^{\theta(\e_i<\e_{i+1})}\over 1-q^2z_i/z_{i+1}}
v_{\cdots ,\e_i,\e_{i+1},\cdots}(\cdots,z_{i+1},z_i,\cdots)\n\\
&&+q{1-z_{i}/z_{i+1}\over 1-q^2z_{i}/z_{i+1}}
v_{\cdots ,\e_{i+1},\e_i,\cdots}(\cdots,z_{i+1},z_i,\cdots),\quad
(\e_i\ne \e_{i+1}),\n\\
&&
v_{\cdots ,\e_i,\e_{i+1},\cdots}(\cdots,z_i,z_{i+1},\cdots)+
q^2{1-q^{-2}z_i/z_{i+1}\over 1-q^2z_{i}/z_{i+1}}
v_{\cdots ,\e_{i+1},\e_i,\cdots}
(\cdots,z_{i+1},z_i,\cdots),\n\\
&&\hskip8cm(\e_i=\e_{i+1}), \label{eq:n}
\end{eqnarray}
where $1\le i\le N-1$. Then  $\N_N=\overline{\hiota_N(\N'_N)}$.

Give $W_N$ a discrete topology.
Set $|0\rangle=v_{\Lambda_0}\otimes 1\in W_0$.
Let $\rho'_N:\V'_N\to W_N$ and  $\rho_N:\V_N\to W_N$
be  the linear maps  defined by 
\begin{eqnarray}
&&v_{\e_1,\cdots,\e_N}(z_1,\cdots,z_N) \mapsto 
\tPhi^{(N-1)}_{\e_1}(z_1)\cdots\tPhi^{(0)}_{\e_N}(z_N)|0\rangle
\,\prod_{0\le k \le N-1} z_{N-k}^{s_k} /\prod_{i<j}\eta(z_j/z_i),\n\\
&&\label{eq:map}
\end{eqnarray}
where $s_k$ is the integer part of $k/(n+1)$ and
\begin{eqnarray}
&&\eta(z)={\displaystyle {(q^2pz;p)_\infty\over (pz;p)_\infty}},\quad
(z;p)_\infty=\prod_{m=0}^\infty(1-p^mz).
\end{eqnarray}
Thanks to the commutation
relations among the intertwiners $\tPhi^{(m)}$, both   maps 
$\rho'_N$ and $\rho_N$ are continuous
and extend to the continuous linear maps 
$\hrho'_N:\hV'_N\to W_N$ and $\hrho_N:\hV_N\to W_N$, respectively.
Moreover the latter maps induce 
$\trho'_N:\hV'_N/\N'_N\to W_N$ and 
$\trho_N:\hV_N/\N_N\to W_N$.

\subsection{Level 0  $U_q'(\widehat{sl_{n+1}})$ actions by Jimbo, et al.} 

Let $U'_q(\widehat{sl_{n+1}})$ be the subalgebra of $U_q(\widehat{sl_{n+1}})$
generated by $e_i$, $f_i$ and  $k_i^{\pm 1}$ $(0\le i\le n)$.
Set $U'_{c=0}=U'_q(\widehat{sl_{n+1}})/<k_0\cdots k_n-1>$ and give this algebra
a discrete topology.
$\hV_N$ is given a topological $U'_{c=0}$ module structure 
by the map
\begin{eqnarray}
&&e_0\mapsto q^{N-1}\sum_{j=1}^N 
\Y_j^{-1} 1^{\otimes j-1}\otimes E_{n+1,1}\otimes
(q^{E_{n+1,n+1}-E_{1,1}})^{\otimes N-j},\n\\
&&f_0\mapsto q^{-(N-1)}\sum_{j=1}^N \Y_j
(q^{E_{1,1}-E_{n+1,n+1}})^{\otimes j-1}\otimes E_{1,n+1}\otimes
1^{\otimes N-j},\n\\
&&k_0\mapsto(q^{E_{n+1,n+1}-E_{1,1}})^{\otimes N},\n\\
&&y\mapsto(\pi\otimes\cdots\otimes\pi)\Delta_0^{op(N-1)}(y),\quad 
y=e_i,f_i,k_i, \quad(1\le i\le n)
\end{eqnarray}
where  $\Y_j^{\pm 1}$ is a continuous linear operator  on $\hV_N$
 defined via the action
of the affine Hecke algebra on Laurent polynomials. See \cite{U_0} for
further details. 
It is further shown that the above $U'_{c=0}$ action on $\hV_N$
induces  the $U'_{c=0}$ action on $\hV_N/\N_N$.

In the case $\mu=\m$, we  replace $\tPhi^{(m)}$ and $W_N$ by $\Phi^{(i)}$ and
$V(\Lambda_j)$ ($i\equiv m$, $j\equiv N$ mod $ n+1$) in (\ref{eq:map}).
We denote the thus obtained maps  by the same letters $\rho'_N$, etc..
Set $\hV=\oplus_{N\ge 1} \hV_N$ and $\H=\oplus_{j=0}^n V(\Lambda_j)$.
Let $\hrho:\hV\to \H$ signify 
the  surjective linear map obtained from $\hrho_N$'s.
In \cite{U_0}, by showing that  $\Ker \hrho$ is $U'_{c=0}$ invariant, 
a $U'_{c=0}$ action is defined on $\H$ 
so that $\hrho$ is a homomorphism of $U'_{c=0}$ modules.

\subsection{Level 0 $U_q'(\widehat{sl_{n+1}})$ actions induced by 
the   toroidal actions}

Let $S_N:\hV'_N\to\hV'_N$ be the  homeomorphic linear map defined by
\begin{equation}
v_{\e_1,\cdots,\e_N}(z_1,\cdots z_N)\mapsto
v_{\e_1,\cdots,\e_N}(z_1,\cdots z_N)/\prod_{i<j}\eta(z_j/z_i).
\end{equation}
Let us define $\sigma_N:\U_\mu\to \End_F(\hV'_N)$ by 
\begin{equation}
\sigma_N(y)=S_N^{-1}\circ (\pi_{a_{N-1}}\hotimes \cdots\hotimes \pi_{a_0})
\Delta^{op}_\I{}^{(N-1)}(y) \circ S_N,\quad y\in \U_\mu.
\end{equation}
Then thanks to Lemma \ref{lem:U},
$\hV'_N$ is a $\U_\mu$ module by the map $\sigma_N$.
(Utilizing  Lemma \ref{lemm:nf}  and (\ref{eq:delta}),
$\N'_N$ is shown to be $\U_\mu$ invariant.)

The $\U_\mu$ modules $\hV'_N$ and $W_N$ are topological 
$U'_{c=0}$ modules 
via the algebra homomorphism  from $U'_q(\widehat{sl_{n+1}})$
to $\U_\mu$ determined by
$e_i\mapsto E_{i,0}$, $f_i\mapsto F_{i,0}$, $k_i\mapsto k_i$.
For these two $U'_{c=0}$ modules, we have the following.
\begin{pro}\label{pro:hrho'}
 The map $\hrho'_N:\hV'_N\to W_N$ is a  homomorphism of $U'_{c=0}$ modules.
\end{pro}
\proof.
Set $\tPhi^{op (m)}=\tPhi^{(m)}\circ \sigma$, 
$(\sigma(v\otimes w)=w\otimes v)$.
The linear map 
$1^{\otimes N-m-1}\otimes \tPhi^{op(m)}:
V_{a_{N-1}/p^{s_{N-1}}} \otimes \cdots \otimes V_{a_m/p^{s_m}}\otimes W_m 
\to
V_{a_{N-1}/p^{s_{N-1}}} \otimes \cdots \otimes V_{a_{m+1}/p^{s_{m+1}}}
\otimes W_{m +1}$
extends to the map between the completions and the latter map	
is   shown to be an intertwiner of $\U_\mu$ modules.
Let $\F^{(N)}=\F^{op (N)}_\I$. From  $(1^{\hotimes N}\hotimes \e)
\F^{(N+1)}=\F^{(N)}$ and (\ref{eq:wt}), we get
$\F^{(N+1)}{}^{\pm 1} (u\otimes |0\rangle)
=(\F^{(N)}{}^{\pm 1}u)\otimes |0\rangle$, $(u\in \hV'_N)$.
Moreover $E_{0,m} |0\rangle=F_{0,m} |0\rangle=0$ $(m\ge 0)$ and 
 $\phi_0^{\pm}(z)|0\rangle=|0\rangle$.
Hence the continuous map 
from $\hV'_N$ to the completion of 
$V_{a_{N-1}} \otimes \cdots \otimes V_{a_0}\otimes W_0$ 
defined  by $u\mapsto S_N(u)\otimes |0\rangle$ is $U'_{c=0}$ linear.
From these we obtain the claim.

\subsection{ Relation between two level 0 $U_q'(\widehat{sl_{n+1}})$ 
actions}

\begin{lem}
Let $X_N$ denote the  span of the vectors 
$v_{\e_1}x^{m_1}\otimes \cdots \otimes v_{\e_N}x^{m_N}$ 
$(\e_1\le \cdots \le \e_N, m_i\in \Z)$. Then 
$X_N+\N'_N/\N'_N$ is dense in $\hV'_N/\N'_N$.
\end{lem}
\proof.
It is sufficient to show $\hV'_N=\overline{X_N+\N'_N}$.
Therefore  the claim follows if we show the following.
\begin{equation}
v_{\e_1}x^{m_1}\otimes \cdots \otimes v_{\e_N}x^{m_N}
\in \overline{X_N+\N'_N}. \label{eq:cl}
\end{equation}
Let $\e_i<\e_{i+1}$ and  $Y_N$ be the the span of the coefficients of 
$v_{\cdots ,\e_i,\e_{i+1},\cdots}(\cdots,z_i,z_{i+1},\cdots)$.
Equation (\ref{eq:n}) implies that the coefficients of
$$
(1-(z_{i}/z_{i+1})^m)v_{\cdots ,\e_{i+1},\e_i,\cdots}
(\cdots, z_{i+1}, z_i,\cdots)
$$
belong to $Y_N+\N'_N$ for any integer $m$. Hence
the coefficients of $v_{\cdots ,\e_{i+1},\e_i,\cdots}
(\cdots, z_{i+1}, z_i,\cdots)$ belong to 
$\overline{Y_N+\N'_N}$. 
Using the last argument  repeatedly, we can show (\ref{eq:cl}).

\begin{lem}\label{lem:eq}
Let $\tiota_N:\hV'_N/\N'_N\to \hV_N/\N_N$
be the map induced by $\hiota_N$.
For $\e_1\le \cdots \le \e_N$  and $y=e_i$, $f_i$,  $k_i^{\pm 1}$
($0\le i\le n$),  the following holds.
\begin{equation}
\tiota_N( y v_{\e_1,\cdots,\e_N}(z_1,\cdots, z_N))=
y\tiota_N(v_{\e_1,\cdots,\e_N}(z_1,\cdots, z_N)).\label{eq:e}
\end{equation}
\end{lem}
\proof.
The above equality clearly holds 
except for $e_0$ and $f_0$. 
The case  $e_0$ is shown in Appendix A and the case $f_0$ is similarly
proven.
\vskip3mm
Since $\tiota_N$ is continuous, the above two lemmas prove the following
proposition.
\begin{pro}\label{pro:psi}
The map $\tiota_N:\hV'_N/\N'_N \to \hV_N/\N_N$
is a  homomorphism of $U'_{c=0}$ modules.
\end{pro}

From Propositions  \ref{pro:hrho'} and \ref{pro:psi},
we can show the following.
\begin{thm}
The map  $\hrho_N:\hV_N\to W_N$ is a homomorphism of $U'_{c=0}$ modules.
\end{thm}
\proof.
Since $\trho_N\circ \tiota_N=\trho'_N$, $\trho_N(yw)=y\trho_N(w)$ 
for $y\in U'_{c=0}$ and
$w\in \Im\tiota_N$. $\Im \tiota_N$ is dense in $\hV_N/\N_N$;
$\hV_N/\N_N$
and $W_N$ are  separated topological $U'_{c=0}$ modules;
$\trho_N$ is continuous. Therefore 
$\trho_N$ is a homomorphism of $U'_{c=0}$ modules.

\vskip3mm
From Proposition \ref{pro:vo} (1), we obtain the following.
\begin{cor}
In the case $\mu=\m$ ($p=\m^2$), the map
$\hrho:\hV\to \H=\oplus_{j=0}^nV(\Lambda_j)$ defined 
in subsection 5.2
 is a surjective homomorphism of $U'_{c=0}$ modules.
\end{cor}

\appendix
\section*{Appendix A.  \hskip5mm   Proof of Lemma 3}
\addtocounter{section}{1}
\setcounter{subsection}{0}
\setcounter{equation}{0}
\renewcommand{\theequation}{\Alph{section}\arabic{equation}}

The following is  immediate. See \cite{U_0} for (3) and (4).
\begin{lemm}\label{lemm:k}
For $I\subset\{1,\cdots N\}$, let $K_I$ be the set of generating series
$$
w(z_1,\cdots z_N)=\sum w_{m_1,\cdots,m_N} z_1^{-m_1}
\cdots z_N^{-m_N}
\mbox{ such that  }
w_{m_1,\cdots,m_N}\in \overline{\V_N[[M]]},
$$
where $M={\rm max}\{m_i\}_{i\in I}$.
For $I\subset\{1,\cdots N\}$ and $k\in\{1,\cdots,N\}$,
let $A_{I,k}$ denote the algebra of
formal power series in the variables $z_i/z_k$ ($ i\in I\setminus\{k\}$).

(1) $A_{I,k} K_I\subset K_{I\setminus\{k\}}$. 
Moreover $A_{I,k}$ acts on $K_I$ when  $k\notin I$.

(2) Let $i,j\in I$, $k\in \{1,\cdots,N\}\setminus I$ and  $w\in K_I$.
 Let further $f(z)$ and $g(z)$
be formal power series in $z$. Then for $(x,y)=(z_i/z_j,z_i/z_k)$
 and $(z_i/z_j,z_j/z_k)$, the following holds.
$$
f(x)(g(y)w)=g(y)(f(x)w)=(f(x)g(y))w.
$$
(3)
Let  $K_{ij}$ be the operator which interchanges $z_i$ and $z_j$,
and set
$\displaystyle \xi_{ij}=q^{-1}+(q-q^{-1}){z_j\over z_i-z_j}(K_{ij}-1)$.
Then  $\xi_{ij}^{\pm 1}K_I\subset K_I$  if $i,j\in I$.

(4) Let $i,j\in I$ and $k\in \{1,\cdots,N\}\setminus I$.
If $f\in A_{I,k}$ satisfies
$f\xi_{ij}=\xi_{ij}f$ when acted on Laurent polynomials
in $z_l$ $(1\le l\le N$),
then $f \xi_{ij}|_{K_I}=\xi_{ij}f|_{K_I}$.

\end{lemm}

\begin{lemm}\label{lemm:l}
Let $L$ be the set of generating series
$$
w(z_1,\cdots z_N)=\sum w_{m_1,\cdots,m_N} z_1^{-m_1}
\cdots z_N^{-m_N}
\mbox{ such that  }
w_{m_1,\cdots,m_N}\in \overline{\V_N[M]},
$$
where $M=\sum_{i=2}^N(i-1)m_i$. 

(1) Let $B$ be the algebra of
formal power series in the variables $z_{i+1}/z_i$ ($1\le i<N$).
Then $B$ acts on $L$.

(2) Set 
$
u_{\e_1,\cdots ,\e_N}(z_1,\cdots, z_N)
=\F^{op (N)}_\I {}^{-1}v_{\e_1,\cdots,\e_N}(z_1,\cdots, z_N).
$

\noindent Then 
$u_{\e_1,\cdots, \e_N}(z_1,\cdots, z_N)\in L$ for any $\e_i$
and 
$u_{\e_1,\cdots ,\e_N}(z_1,\cdots, z_N)=v_{\e_1,\cdots,\e_N}(z_1,\cdots, z_N)$
when $\e_1\le \cdots \le \e_N$.
\end{lemm}
\proof.
(1) is immediate. (2)  follow from (\ref{eq:wt}).

\begin{lemm}\label{lemm:nf}
Let $\N'_N{}^\F$ be the closure in $\hV'_N$ 
of the span of the coefficients of the following generating series,
\begin{eqnarray}
&&v_{\cdots ,\e_i,\e_{i+1},\cdots}
(\cdots,z_i,z_{i+1},\cdots)+
q{1-z_{i}/z_{i+1}\over 1-q^2z_{i}/z_{i+1}}
v_{\cdots ,\e_{i+1},\e_i,\cdots}
(\cdots,z_{i+1},z_i,\cdots),\n\\
&&\hskip10cm (\e_i<\e_{i+1}),\n\\
&&v_{\cdots ,\e_i,\e_{i+1},\cdots}(\cdots,z_i,z_{i+1},\cdots)
+q^2{1-q^{-2}z_i/z_{i+1}\over 1-q^2z_{i}/z_{i+1}}
v_{\cdots ,\e_{i+1},\e_i,\cdots}(\cdots,z_{i+1},z_i,\cdots),\n\\
&&\hskip9cm(\e_i=\e_{i+1}), \label{eq:nf}
\end{eqnarray}
where $1\le i\le N-1$.
Then $\F^{op(N)}\N'_N=\N'_N{}^\F$.
\end{lemm}

\proof.
Let $\F^{(N)}$ signify $\F_\I^{op (N)}$.
Firstly we consider the case $N=2$. Let $\e_1< \e_2$.
From the definition we get 
\begin{equation}
v_{\e_1,\e_2}(z_1,z_2)+q{1-z_2/z_1\over 1-q^2z_2/z_1} v_{\e_2,\e_1}(z_1,z_2)
+(z_1\leftrightarrow z_2)\in \N'_N{}^\F.
\end{equation}
This implies
\begin{equation}
{z_2\over 1-z_2/z_1}v_{\e_1,\e_2}(z_1,z_2)+{qz_2\over 1-q^2z_2/z_1}
v_{\e_2,\e_1}(z_1,z_2)-(z_1\leftrightarrow z_2)\in \N'_N{}^\F.
\end{equation}
From this and the  equality
\begin{eqnarray}
\F^{(2)} v_{\e_1,\e_2}(z_1,z_2)
&=&v_{\e_1,\e_2}(z_1,z_2),\quad (\e_1\le\e_2),\n\\
&=&v_{\e_1,\e_2}(z_1,z_2)-
(q-q^{-1}){z_2/z_1\over 1-z_2/z_1}v_{\e_2,\e_1}(z_1,z_2),\quad (\e_1>\e_2),\n\\
&&
\end{eqnarray}
we can show that the case $N=2$ holds. 
The case $N>2$ follows from the fact that 
$\N'_2$ is $U$ invariant and the  equality
\begin{equation}
\F^{(N)}= \F^{(2)}_{k,k+1} 
(1^{\hotimes k-1}\hotimes \Delta_\I^{op}\hotimes 1^{\hotimes N-k-1})
\F^{(N-1)}.
\end{equation}

Hereafter we shall let $\equiv$ denote the equality in $\hV_N$ 
mod $\N_N$.
\begin{lemm}\label{lemm:rep}
The action of $e_0$ on the $U'_{c=0}$ module $\hV_N$ 
and that of $E_{0,0}$  on the $\U_\mu$ module $\hV'_N$
satisfy the following relations.

\begin{eqnarray}
(1)&&e_0 v_{\e_1,\cdots,\e_N}(z_1,\cdots,z_N)\n\\
&\equiv&(-q)^{N-1}\sum_{j=1}^N \delta_{\e_j,1}(-1)^{j-1}
q^{-\sum_{j<i\le N}\delta_{\e_i,1}}
\xi_{j-1\,j}\cdots\xi_{1j}\n\\
&&\times v_{\e_1,\cdots,\hat{\e_j},\cdots, \e_N,n+1}
(z_1,\cdots,\hat{z_j},\cdots, z_N, z_j/p). \\
(2)&&\mbox{ For  }\e_1\le \cdots \le\e_N,\n\\
&&\iota_N(E_{0,0}v_{\e_1,\cdots,\e_N}(z_1,\cdots, z_N))\n\\
&\equiv&(-q)^{N-1}\sum_{j=1}^N\delta_{\e_j,1}(-q^{-1})^{j-1}
 q^{-\sum_{j<i\le N} \delta_{\e_i,1}}
\prod_{i=1}^{j-1}{1-q^2z_j/z_i\over 1-z_j/z_i}
\prod_{i=j+1}^N({1-q^2z_i/z_j\over 1-z_i/z_j})^{\delta_{\e_i,1}}\n\\
&&\times 
v_{\e_1,\cdots,\hat{\e_j},\cdots,\e_N,n+1}(z_1,\cdots,\hat{z_j},
\cdots, z_N, z_j/p).
\end{eqnarray}
Here $\,\hat{}\,$ denotes the omission of variables.
\end{lemm}

\proof.
(1) follows from the definition of $\Y_j^{-1}$ \cite{U_0}, (\ref{eq:n})
and Lemma \ref{lemm:k}.
(2) Set $\displaystyle f(z)=q^{-1}{1-q^2z\over 1-z}$. Set further
$\displaystyle g_\e(z)=q{1-z\over 1-q^2z}$ for $1\le \e\le n$ and 
$\displaystyle =q^2{1-q^{-2}z\over 1-q^2z}$ for $\e=n+1$.
Then thanks to  (\ref{eq:delta}) we get
\begin{eqnarray}
E_{0,0}u_{\e_1,\cdots \e_N}(z_1,\cdots, z_N)&=&q^{N-1}
\sum_{j=1}^N \prod_{i=1}^{j-1}f(z_j/z_i)
\prod_{i=j+1}^{N}f(z_i/z_j)^{\delta_{\e_i,1}}\n\\
&\times&
\prod_{i=j+1}^{N}g_{\e_i}(pz_i/z_j)
u_{\e_1,\cdots,\e_{j-1},n+1,\e_{j+1},\cdots, \e_N}
(z_1,\cdots, z_j/p,\cdots z_N),\n\\
&&
\end{eqnarray}
for any $\e_i$. Since $\iota_N$ is continuous, 
Lemmas \ref{lemm:l} and  \ref{lemm:nf} imply
\begin{eqnarray}
&&\prod_{i=j+1}^{N}g_{\e_i}(pz_i/z_j)
\iota_N(u_{\e_1,\cdots,\e_{j-1},n+1,\e_{j+1},\cdots, \e_N}
(z_1,\cdots, z_j/p,\cdots z_N))\n\\
&\equiv&(-1)^{N-j}\iota_N(u_{\e_1,\cdots,\hat{\e_j},\cdots, \e_N, n+1}
(z_1,\cdots, \hat{z_j}, \cdots z_N, z_j/p)).\n\\
&&
\end{eqnarray}
In the case $\e_1\le \cdots \le \e_N$, 
thanks to Lemma \ref{lemm:l} (2), multiplying the above equation 
by  $\prod_{i=1}^{j-1}f(z_j/z_i)
\prod_{i=j+1}^{N}f(z_i/z_j)^{\delta_{\e_i,1}}$ is meaningful.
Hence we obtain the claim.

\vskip3mm

Therefore  (\ref{eq:e}) with $y=e_0$
follows from the following lemma with 
$w(z_1,\cdots, z_N)=
v_{\underbrace{\scriptstyle{1,\cdots, 1}}_{s-1},\e_{s+1},
\cdots,\e_N, n+1}(z_1,\cdots z_N)$.
\begin{lemm}\label{lemm:a}
Let $1\le t\le s\le N$ and $w(z_1,\cdots, z_N)\in 
K_{\{t, \cdots ,N\}}$.
If the coefficients of the  generating series
$$
w(z_1,\cdots,z_i,z_{i+1},\cdots,z_N)
+q^2{1-q^{-2}z_{i}/z_{i+1}\over 1-q^2z_{i}/z_{i+1}}
w(z_1,\cdots,z_{i+1},z_i,\cdots, z_N),\quad (t\le i\le s-2)
$$
belong to $\N_N$,
then
the following equality holds   in $\hV_N$.
\begin{eqnarray}
&&\sum_{j=t}^s(-q)^{j-t} \xi_{j-1\,j}\cdots \xi_{tj}
w(z_1,\cdots,\hat{z_j},\cdots,z_N,z_j/p)\n\\
&\equiv&\sum_{j=t}^s(-1)^{j-t}\prod_{i=t}^{j-1}{1-q^2z_j/z_i\over 1-z_j/z_i}
\prod_{i=j+1}^s{1-q^2z_i/z_j\over 1-z_i/z_j}
w(z_1,\cdots,\hat{z_j},\cdots,z_N,z_j/p).\n\\
&&
\end{eqnarray}
\end{lemm}

\proof. 
This can be shown by induction on $s-t$.
Thanks to Lemma \ref{lemm:k}, the left hand side is rewritten as follows,
\begin{eqnarray}
&&w(z_1,\cdots, \hat{z_t},\cdots, z_N,z_t/p)
+(1-q^2)\sum_{j=t+1}^s(-q)^{j-t-1}\xi_{j-1\,j}\cdots\xi_{t+1\,j}\n\\
&&\times{z_j/z_t\over 1-z_j/z_t}
w(z_1,\cdots,z_{t-1},z_j,z_{t+1},\cdots, \hat{z_j},\cdots,z_N,z_t/p)\n\\
&&-{1\over \prod_{i>t}\eta(z_i/z_t)}
\sum_{j=t+1}^s(-q)^{j-t-1}\xi_{j-1\,j}\cdots\xi_{t+1\,j}
 {\bar w}(z_1,\cdots, \hat{z_j},\cdots,z_N,z_j/p),\n\\
&&
\end{eqnarray}
where
\begin{equation}
{\bar w}(z_1,\cdots, z_N)=\prod_{i>t}\eta(z_i/z_t)\times w(z_1,\cdots,z_N)
\in K_{\{t+1,\cdots, N\}}.
\end{equation}
For $t+1\le i<j\le s$,  thanks to Lemma \ref{lemm:k} (2), we obtain,
\begin{eqnarray}
&&\xi_{ij}{ z_j/z_t\over 1- z_j/z_t}
w(z_1,\cdots,\hat{z_t},\cdots,
z_{i-1}, z_j,z_i, \cdots,\hat{z_j}, \cdots,z_N, z_t/p)\n\\
&\equiv&-q^{-1}{1- q^2z_i/z_t\over 1-z_i/z_t}{ z_j/z_t\over 1- z_j/z_t}
w(z_1,\cdots,\hat{z_t},\cdots,
z_i,z_j,z_{i+1}, \cdots,\hat{z_j}, \cdots, z_N, z_t/p).\n\\
&&
\end{eqnarray}
Since the action of $\xi_{j-1\,j} \cdots \xi_{i+1\,j}$ on
both hand sides of the above is well defined, we get the following equality.
\begin{eqnarray}
&&\xi_{j-1\,j}\cdots \xi_{ij}{z_j/z_t\over 1- z_j/z_t}
w(z_1,\cdots, \hat{z_t}, \cdots, z_{i-1}, z_j,z_i,
\cdots,\hat{z_j},\cdots,  z_N, z_t/p)\n\\
&\equiv&-q^{-1}{1- q^2z_i/z_t\over 1-z_i/z_t}
\xi_{j-1\,j}\cdots \xi_{i+1\, j}{z_j/z_t\over 1- z_j/z_t}
w(z_1,\cdots, \hat{z_t},\cdots, z_i,z_j,z_{i+1},
 \cdots,\hat{z_j},\cdots,  z_N, z_t/p).\n\\
&&
\end{eqnarray}
Repeating this argument, the sum of the first two terms is 
found to be
$$
\prod_{i=t+1}^s{1- q^2z_i/z_t\over 1- z_i/z_t}
w(z_1,\cdots,\hat{z_t},\cdots, z_N, z_t/p).
$$
Applying the assumption of the induction 
to the last sum, we obtain the claim.

\end{document}